\documentclass{colt2012} 

\usepackage{algorithm,algorithmic}
\usepackage{amssymb, multirow, paralist, color}
\newtheorem{thm}{Theorem}

\title[Unified  Convergence Analysis of Stochastic Momentum Methods ]{Unified  Convergence Analysis of Stochastic Momentum Methods for Convex and Non-convex Optimization}
 \coltauthor{\Name{Tianbao Yang}$^\dagger$ \Email{tianbao-yang@uiowa.edu}\\
 \Name{Qihang Lin}$^\ddagger$ \Email{qihang-lin@uiowa.edu}\\
\Name{Zhe Li}$^\dagger$\Email{zhe-li-1@uiowa.edu}\\
   \addr$^\dagger$Department of Computer Science \\
    The University of Iowa, Iowa City, IA 52242 \\
   \addr$^\ddagger$Department of Management  Sciences \\
    The University of Iowa, Iowa City, IA 52242
}

\def \y {\mathbf{y}}
\def \E {\mathrm{E}}
\def \x {\mathbf{x}}

\def \z {\mathbf{z}}

\def \R {\mathbb{R}}

\def \v {\mathbf{v}}

\def \p {\mathbf{p}}

\def \xh {\widehat{\x}}

\def \G {\mathcal G}

\begin{document}
\maketitle

\begin{abstract}
Recently, {\it stochastic momentum} methods have been  widely adopted in training deep neural networks. However, their convergence analysis is  still underexplored at the moment, in particular for non-convex optimization. This paper fills  the gap between practice and theory by developing a basic convergence analysis of two stochastic momentum methods, namely stochastic heavy-ball method and the stochastic variant of Nesterov's accelerated gradient method. We hope that the basic convergence results developed in this paper can serve the reference to the convergence of stochastic momentum methods and also serve the baselines for comparison in future development of stochastic momentum methods. The novelty of convergence analysis presented in this paper is a unified framework, revealing  more insights about the similarities and differences between different stochastic momentum methods and stochastic gradient method. The unified framework exhibits a continuous change  from the gradient method to Nesterov's accelerated gradient method and finally the heavy-ball method incurred by a free parameter, which can help explain a similar change observed  in the testing error convergence behavior for deep learning. Furthermore, our empirical results  for optimizing deep neural networks demonstrate that the stochastic variant  of Nesterov's accelerated gradient method achieves a good tradeoff  (between speed of convergence  in training error and  robustness of convergence in testing error) among the three stochastic methods.

\end{abstract}

\section{Introduction}
Momentum methods have a long history dating back to 1960's. \cite{poly64} proposed a heavy-ball (HB) method that uses the previous two iterates when computing the next one. For a twice continuously differentiable strongly convex and smooth objective function, Polyak's analysis yields an accelerated  linear convergence rate. In 1983, \cite{citeulike:9501961} proposed an accelerated gradient (NAG) method, which achieves the optimal $O(1/t^2)$ convergence rate for convex smooth optimization~\footnote{$t$ is the number of iterations.}. It was later on shown to have an accelerated linear convergence rate for smooth and strongly convex optimization problems~\citep{opac-b1104789}. Both the HB method and the NAG method use a momentum term in updating the solution, i.e., the difference between current iterate and the previous iterate. Therefore, both methods have been referred to as momentum methods in literature.  

Although the convergence analysis of HB has been recently extended to smooth functions for both convex~\citep{arxiv1412,DBLP:journals/jmiv/OchsBP15} and non-convex optimization~\citep{DBLP:journals/siamis/OchsCBP14}, it is still under-explored for stochastic optimization, specially stochastic non-convex optimization. In practice, the stochastic variants of HB and NAG that use a (mini-batch) stochastic gradient or subgradient  have been employed widely in deep learning, leading to tremendous success in many areas~\citep{krizhevsky2012imagenet}.  Due to that training a  deep neural networks (DNN) usually requires a huge amount of data,  the stochastic variants are much more efficient than their deterministic couterparts. Recent studies also find that momentum term is useful in deep network training and more generally non-convex optimization~\citep{DBLP:conf/icml/SutskeverMDH13,DBLP:journals/siamis/OchsCBP14}. The illustrative example studied  in~\citep{DBLP:journals/siamis/OchsCBP14} shows that the momentum term has the potential to escape from the local minimum in the deterministic setting. 

Nonetheless,  there still lacks any convergence analysis of the stochastic momentum methods for convex and non-convex optimization.   Although one may believe that they should entail similar convergence rate as the stochastic gradient (SG) methods, it is still of interest to develop their own convergence analysis for  convex and non-convex optimization. To bridge this gap between  practice and theory, we develop a basic convergence analysis of both the stochastic HB method and the stochastic NAG method in a unified framework. Without any surprise, for a general Lipschitz continuous convex objective function we show that the objective value converges in expectation at the order of $1/\sqrt{t}$, where $t$ is the number of iterations, and for a non-convex objective function  with Lipschitz continuous gradient we show that the gradient converges in expectation at the order of $1/\sqrt{t}$. The novelty of convergence analysis presented in this paper is a unified framework, revealing  more insights about the similarities and differences between different stochastic momentum methods and stochastic gradient method. In particular, the unified framework exhibits a continuous change  from the gradient method to Nesterov's accelerated gradient method and finally the heavy-ball method incurred by a free parameter, which can help explain a similar change observed  in the testing error convergence behavior for deep learning. 

It is worth mentioning that several recent works by \cite{DBLP:journals/mp/GhadimiL16,DBLP:journals/siamjo/GhadimiL13a} have established similar  convergence results of the stochastic gradient method and the stochastic version of a variant of accelerated gradient method for non-convex optimization under the same setting. However, the variant of accelerated gradient method in~\citep{DBLP:journals/mp/GhadimiL16} is hard to be explained in the framework of momentum methods and their analysis is not applicable to the stochastic heavy-ball method.  Additionally, two recent studies~\citep{DBLP:journals/corr/ReddiHSPS16,DBLP:journals/corr/ZhuH16a} have analyzed stochastic variance reduced gradient (SVRG) methods for nonconvex finite-sum problems  and established potentially faster convergence than SG. It might be interesting to analyze  stochastic momentum methods using variance reduced stochastic gradients by exploring  our analysis and the analysis of SVRG for non-convex finite-sum problems, which is beyond the scope of this paper and is left as a future work. 


In the reminder of the paper, we will first review the HB method and the NAG method and present their stochastic variants. Then we present a unified view of these momentum methods. In Section~\ref{sec:conv}, we present the convergence analysis for stochastic momentum methods. In Section~\ref{sec:exp}, we present some empirical results for comparing different methods for optimizing deep neural networks. Finally, we conclude this work. 

\section{Momentum Methods and their Stochastic variants}
The problem of interest  can be cast into 
\begin{align}
\min_{\x\in\R^d} f(\x)
\end{align}
where $f(\x)$ is not necessarily a convex function. We use $\nabla f(\x)$ to denote the gradient of a smooth function. A function is smooth iff there exists $L>0$ such that
\begin{align}
\|\nabla f(\y) - \nabla f(\x)\|\leq L\|\y - \x\|, \quad \forall \x, \y\in\R^d
\end{align}
where $\|\cdot\|$ denotes the Euclidean norm. Note that the above inequality does not imply the convexity.  For a non-smooth function $f(\x)$, we denote by $\partial f(\x)$ the set of its subgradients. A function $f(\x)$ is convex iff there exists $\v\in\partial f(\x)$ such that , 
\begin{align}
f(\y)\geq f(\x) + \v^{\top}(\y - \x)
\end{align}

There are two types of momentum methods, one is based on the HB method and the other is based on Nesterov's  AG (NAG) method.  The HB method was originally proposed for optimizing a smooth and strongly convex objective function. Its update is given below for $k=0,\ldots, $
\begin{align}\label{eqn:hb}
\text{HB:} \quad \x_{k+1} = \x_k - \alpha\nabla f(\x_k) + \beta(\x_k - \x_{k-1})
\end{align}
with $\x_{-1} = \x_0$, where $\beta\in[0,1)$ is called the momentum constant and $\alpha$ is the step size. Equivalently, the above update can be implemented by the following two steps for $k=0,\ldots, $: 
\begin{equation}\label{eqn:hb2}
\hspace*{-0.3in}\text{HB:}\quad\left\{ \begin{aligned}
\v_{k+1} &= \beta \v_{k} - \alpha \nabla f(\x_k)\\
\x_{k+1} &  =  \x_k+\v_{k+1}
\end{aligned}\right.
\end{equation}
The update of NAG consists of the following two steps for $k=0,\ldots,$
\begin{equation}\label{eqn:ag1}
\text{NAG:}\quad \left\{\begin{aligned}
\y_{k+1} & = \x_k - \alpha \nabla f(\x_k)\\
\x_{k+1}& = \y_{k+1} + \beta(\y_{k+1} - \y_k)
\end{aligned}\right.
\end{equation}
with $\y_0 = \x_0$. 
By introducing $\v_k = \y_k - \y_{k-1}$ with $\v_0 =0$, the above update can be equivalently written as 
\begin{equation}\label{eqn:ag}
\text{NAG:}\quad \left\{\begin{aligned}
 \v_{k+1}& = \beta\v_k - \alpha \nabla f(\y_k + \beta \v_k)\\
\y_{k+1}&  =\y_k + \v_{k+1}
\end{aligned}\right.
\end{equation}
By comparing~(\ref{eqn:ag}) to~(\ref{eqn:hb2}), one might argue that the difference between HB and NAG lies at the point for evaluating the gradient~\citep{DBLP:conf/icml/SutskeverMDH13}. We will present a unified view of the two methods that allows us to analyze them in a unified framework. It is notable that if the momentum constant $\beta=0$, both HB and NAG reduce to the gradient method. The convergence of HB and NAG has been established for convex optimization~\citep{poly64,citeulike:9501961,opac-b1104789,arxiv1412,DBLP:journals/jmiv/OchsBP15}. One issue of the updates in HB and NAG  is that they need to compute the gradient of the objective function, which in some applications (e.g., machine learning) is prohibitive  for big data. Therefore, when employed for optimizing deep neural networks, the gradient is usually replaced with the stochastic (sub)-gradient, which yields the stochastic variants. 
We denote by $\G_k=\G(\x_k;\xi_k)$ a stochastic gradient (resp. subgradient) of $f(\x)$ at $\x_k$ depending on a random variable $\xi_k$ such that  $\E[\G(\x_k ;\xi_k)]=\nabla f(\x_k)$ (resp. $\E[\G(\x_k ;\xi_k)]\in\partial f(\x_k)$). Then the update of stochastic HB (SHB) becomes 
\begin{align}\label{eqn:shb}
\text{SHB:} \quad \x_{k+1} = \x_k - \alpha \G(\x_k; \xi_k) + \beta(\x_k - \x_{k-1})
\end{align}
The update of stochastic NAG (SNAG) becomes 
\begin{equation}\label{eqn:sag}
\hspace*{-0.65in}\text{SNAG:}\quad \left\{\begin{aligned}
 \y_{k+1}& = \x_k - \alpha \G(\x_k; \xi_k)\\
\x_{k+1}&  =\y_{k+1} + \beta(\y_{k+1} - \y_k)
\end{aligned}\right.
\end{equation}

\subsection{A Unified View of Stochastic Momentum Methods}
In this subsection, we present a unified view of the two (stochastic)  momentum methods and  (stochastic) (sub)-gradient methods. Since the choice of (sub)-gradient and stochastic (sub)-gradient is irrelevant to our discussion here, we denote by $\G(\x)$ either a (sub)-gradient or a stochastic (sub)-gradient. We first present the unified framework and then show that HB, NAG and the (sub)-gradient method are special cases of the unified framework. Let $\alpha>0$, $\beta\in[0,1)$, and  $s\geq 0$. The updates of the unified momentum (UM) method  are given by 
\begin{equation}\label{eqn:um}
\text{UM}:\quad\left\{\begin{aligned}
\y_{k+1} &  = \x_k - \alpha \G(\x_k)\\
\y^s_{k+1} & = \x_k - s\alpha \G(\x_k)\\
\x_{k+1} & = \y_{k+1} + \beta(\y^s_{k+1} - \y^s_k)
\end{aligned}\right.
\end{equation}
with $\y^s_0 = \x_0$. We denote by SUM the stochastic version of UM when using a stochastic (sub)-gradient in the above updates.  It is notable that in the update of $\x_{k+1}$, a momentum term is constructed based the auxiliary sequence  $\{\y_k^s\}$, which are intermediate updates based on $\x_k$  using a possibly different step size from that for updating $\y_{k}$. In particular, HB uses a zero step size to construct the auxiliary sequence $\y^s_{k}$, NAG uses a step size $\alpha$ the same to that for $\y_k$ and the gradient method uses a larger step size $\frac{\alpha}{1-\beta}$. 

Next, we discuss three choices for $s=0, 1$ and $s = \frac{1}{1-\beta}$, which correspond to HB, NAG and (sub)-gradient method, respectively. 

When $s=0$, then $\y^s_{k+1} = \x_k$, the update of $\x_{k+1}$ is equivalent to 
\[
\x_{k+1} = \x_k - \alpha \G(\x_k) + \beta(\x_{k} - \x_{k-1})
\]
which is exactly the update of the HB method. 

When $s=1$, then $\y^s_{k+1} = \y_{k+1}$, then the update in~(\ref{eqn:um}) reduces to that in~(\ref{eqn:ag1}) or~(\ref{eqn:sag}) of the AG method. 

The last special case corresponds to using $s = \frac{1}{1-\beta}$. We will show that the update of $\x_{k+1}$ is equivalent to 
\begin{align}\label{eqn:gd}
\x_{k+1} = \x_k - \frac{\alpha}{1 - \beta}\G(\x_k), k\geq 0
\end{align}
which is exactly the update of the gradient method using a step size $\frac{\alpha}{1 - \beta}$. First, we verify~(\ref{eqn:gd}) holds for $k=0$. From the updates in~(\ref{eqn:um}), we have
\begin{align*}
\x_1  & = \y_1 + \beta(\y^s_1 - \y^s_0) = \x_0 - \alpha \G(\x_0) + \beta(\x_0 - s\alpha \G(\x_0) - \x_0)\\
& =  \x_0 - \alpha \G(\x_0) - s \beta \alpha \G(\x_0) = \x_0 - \alpha \G(\x_0)(1 + \frac{\beta}{1-\beta}) = \x_0 - \frac{\alpha}{1-\beta}\G(\x_0)
\end{align*}
Then we show~(\ref{eqn:gd}) holds for any $k\geq 1$. From  the updates in~(\ref{eqn:um}), we have
\begin{align*}
\x_{k+1} - \x_k = - \alpha \G(\x_k) + \beta (\y^s_{k+1} - \y^s_k)  = - \alpha \G(\x_k) + \beta(\x_k - s\alpha\G(\x_k) - \x_{k-1} + s\alpha \G(\x_{k-1}))
\end{align*}
Then
\begin{align*}
\x_{k+1} - \x_k + s\alpha  \G(\x_k) = \beta(\x_k - \x_{k-1} + s\alpha \G(\x_{k-1})) + (s - 1 - \beta s)\alpha \G(\x_k)
\end{align*}
Since $s = \frac{1}{1-\beta}$, then $s - 1 - \beta s=0$, thus
\begin{align*}
\x_{k+1} - \x_k + s\alpha  \G(\x_k) = \beta(\x_k - \x_{k-1} + s\alpha \G(\x_{k-1})), k\geq 1
\end{align*}
Therefore
\[
\x_{k+1} - \x_k + s\alpha \G(\x_k)= \beta^{k} (\x_1 - \x_0 + \frac{\alpha}{1-\beta}\G(\x_0) ) = 0 
\]
which leads to~(\ref{eqn:gd}). 

To facilitate the unified analysis of the stochastic momentum methods, we note that~(\ref{eqn:um}) implies the following recursions, which are straightforward to verify 
\begin{align}\label{eqn:rec}
\x_{k+1} + \p_{k+1} & = \x_k + \p_k  - \frac{\alpha}{ 1- \beta}\G(\x_k), \quad k\geq 0\\
\v_{k+1} & = \beta \v_k + ((1-\beta)s - 1) \alpha \G(\x_k),\quad k\geq 0\label{eqn:rec2}
\end{align}
where $\v_k = \frac{(1-\beta)}{\beta}\p_k$ and $\p_k$ is  given by 
\begin{equation}\label{eqn:p}
\p_k = \left\{
\begin{aligned}
&\frac{\beta}{1-\beta}(\x_k - \x_{k-1} + s\alpha \G(\x_{k-1})),\quad k\geq 1\\
& 0,\quad k=0\\
\end{aligned}\right.
\end{equation}
It is worth mentioning that the recursion in~(\ref{eqn:rec}) has been observed and employed  in~\citep{arxiv1412}, but it focused on the convergence analysis of the heavy-ball method for convex optimization in the deterministic setting. 


\section{Unified Convergence Analysis for Convex and Non-convex Optimization}\label{sec:conv}
In this section, we present the basic convergence results and analysis of the stochastic momentum methods. 
We first state the main results and then present the analysis.  In the sequel, we let $\G_k = \G(\x_k;\xi_k)$ denote a stochastic gradient or a stochastic subgradient. 

\subsection{Main results}
\subsubsection{Convex Optimization}
\begin{thm}~\label{thm:1} (Convergence of SUM) 
Suppose $f(\x)$ is a convex function, $\E[\|\G(\x; \xi) - \E[\G(\x; \xi)]\|^2]\leq \delta^2$ and $\|\partial f(\x)\|\leq G$ for any $\x$. 
   Let update~(\ref{eqn:um}) run for $t$ iterations with $\G(\x_k;\xi_k)$. By setting $\alpha=\frac{C}{\sqrt{t+1}}$ we have
 \[
 \E[f(\xh_t) - f(\x_*)]\leq \frac{\beta (f(\x_0) - f(\x_*)) }{(1-\beta)(t+1)}+  \frac{(1-\beta)\|\x_0 - \x_*\|^2}{2C\sqrt{t+1}} + \frac{C(1+2s\beta)(G^2 +\delta^2)}{2(1-\beta)\sqrt{t+1}}
 \]
 where $C$ is a postive constant, $\xh_{t} = \sum_{k=0}^t\x_k/(t+1)$ and $\x_*$ is an optimal solution.
\end{thm}
{\bf Remark:} The difference of the convergence bounds for different variants of SUM only lies at the constant factor $1+2s\beta$ in the third term.  We emphasize here that this constant factor $1+2s\beta$  is artifact of our analysis, which should not be understood that SNAG or SG is worse than SHB. Indeed, with a refined analysis, we might achieve a better convergence for SNAG (more generally SUM when $s\beta\geq 1/2$) under a slightly different  condition.

 \begin{thm}~\label{thm:2}
Suppose $f(\x)$ is a convex function and $\|\G(\x; \xi) \|\leq M$.   Let update~(\ref{eqn:um}) run for $t$ iterations with $s\beta\geq 1/2$. We  have
\begin{align*}
\E[f(\xh_t) - f(\x)]&\leq \frac{\beta(f(\x_0) - f(\x_*))}{(1-\beta)(t+1)} +\frac{(1-\beta) \|\x_0 - \x_*\|^2}{2\alpha(t+1)} + \frac{s\alpha\beta M}{(1-\beta)}\frac{\E\left[\sum_{k=0}^t\|\G_k - \G_{k-1}\|\right]}{t+1}
\end{align*}
 where $C$ is a postive constant, $\xh_{t} = \sum_{k=0}^t\x_k/(t+1)$ and $\x_*$ is an optimal solution.
\end{thm}
{\bf Remark:} The third part in the above convergence bound has a variational term in terms of the stochastic subgradient  $V_t=\E\left[\sum_{k=0}^t\|\G_k - \G_{k-1}\|\right]$. In the worse case, $V_t\leq 2M(t+1)$ which yields a bound with a similar order to that in Theorem~\ref{thm:1} by setting $\alpha = \frac{C}{\sqrt{t+1}}$. If $V_t\leq o(t)$, we could achieve a better convergence.

\subsubsection{Non-convex Optimization}

\begin{thm}~\label{thm:3}(Convergence of SUM) 
Suppose $f(\x)$ is a non-convex and $L$-smooth function, $\E[\|\G(\x; \xi) - \nabla f(\x)\|^2]\leq \delta^2$ and $\|\nabla f(\x)\|\leq G$ for any $\x$. 
 Let update~(\ref{eqn:um}) run for $t$ iterations with $\G(\x_k; \xi_k)$. By setting $\alpha =\min\{\frac{1-\beta}{2L}, \frac{C}{\sqrt{t+1}}\}$ we have
\begin{align*}
\min_{k=0,\ldots, t}\E[\|\nabla f(\x_k)\|^2]&\leq \frac{2(f(\x_0) - f_*)(1-\beta)}{t+1}\max\left\{\frac{2L}{1-\beta}, \frac{\sqrt{t+1}}{C}\right\} \\
&+ \frac{C}{\sqrt{t+1}}\frac{L\beta^2((1-\beta)s -1)^2(G^2+\sigma^2) + L\sigma^2(1-\beta)^2}{(1-\beta)^3} \end{align*}
\end{thm}
{\bf Remark:}  We would like to make several remarks. First, the difference of the convergence bounds for different variants of SUM lies at term $L\beta^2((1-\beta)s-1)^2(G^2+\sigma^2))$, which is equal to $L\beta^2$, $L\beta^4$ and $0$ for SHB, SNAG and SG, respectively. Second,  two different views of SG yield the same convergence result as that in~\citep{DBLP:journals/siamjo/GhadimiL13a}. The first view considers SG as a special  case of of SUM with $\beta=0$, which gives the 
update $\x_{k+1} = \x_{k} - \alpha \G_k$. Theorem~\ref{thm:3} then implies a convergence bound of 
\begin{align}\label{eqn:bsg}
\frac{2(f(\x_0) - f_*)}{t+1}\max\left\{2L, \frac{\sqrt{t+1}}{C}\right\} + \frac{CL\sigma^2}{\sqrt{t+1}}
\end{align}
with $\alpha = \min\{1/(2L), C/\sqrt{t+1}\}$.  The second view considers SG as a special case of SUM with $s = \frac{1}{1-\beta}$, which yields the update in~(\ref{eqn:gd}), i.e., $\x_{k+1} = \x_k - \frac{\alpha}{1-\beta}\G_k$. Theorem~\ref{thm:3} then implies the same convergence bound in~(\ref{eqn:bsg}) with $\alpha = \min\{\frac{1-\beta}{2L}, \frac{C(1-\beta)}{\sqrt{t+1}}\}$. Third, the step size $\alpha$ of different variants of SUM used in the analysis of Theorem~\ref{thm:3} is the same value. Below, we present a result with different step sizes $\alpha$ for different variants of SUM in the analysis, which sheds more insights of different methods. 
\begin{thm}~\label{thm:4}(Convergence of SUM) 
Suppose $f(\x)$ is a non-convex and $L$-smooth function, $\E[\|\G(\x; \xi) - \nabla f(\x)\|^2]\leq \delta^2$ and $\|\nabla f(\x)\|\leq G$ for any $\x$. 
Let update~(\ref{eqn:um}) run for $t$ iterations with $\G(\x_k; \xi_k)$. By setting $\alpha =\min\{\frac{1-\beta}{2L[1+((1-\beta)s-1)^2]}, \frac{C}{\sqrt{t+1}}\}$ we have
\begin{align*}
\min_{k=0,\ldots, t}\E[\|\nabla f(\x_k)\|^2]&\leq \frac{2(f(\x_0) - f_*)(1-\beta)}{t+1}\max\left\{\frac{2L[1+((1-\beta)s-1)^2]}{1-\beta}, \frac{\sqrt{t+1}}{C}\right\} \\
&+ \frac{C}{\sqrt{t+1}}\frac{L\beta^2(G^2+\sigma^2) + L\sigma^2(1-\beta)^2}{(1-\beta)^3} \end{align*}
\end{thm}
{\bf Remark: } The above result allows us to possibly set a larger initial  step size for SG ($s=1/(1-\beta)$) and SNAG ($s=1$) than that for SHB ($s=0$). Our empirical studies for deep learning in Section~\ref{sec:exp} also demonstrate this point. 


\subsection{Analysis}
\subsubsection{Proof of Thereom~\ref{thm:1}}

\begin{lemma}\label{lem:0}
Let $\x_{-1}=\x_0$. For any $k\geq 0$, we have
\begin{align*}
\E[\|\x_{k+1} + \p_{k+1} - \x\|^2&] \leq   \E[\|\x_k + \p_k - \x\|^2] - \frac{2\alpha}{1-\beta}\E[(f(\x_k) - f(\x))]\\
 &- \frac{2\alpha\beta}{(1-\beta)^2}\E[(f(\x_k) - f(\x_{k-1}))]  + \left(\frac{\alpha}{1-\beta}\right)^2(2s\beta+1) (G^2 +\delta^2)
\end{align*}
\end{lemma}
Summing the above inequality for $k=0, \ldots, t$ gives 
\begin{align*}
\frac{2\alpha}{1-\beta}\sum_{k=0}^t\E[f(\x_k) - f(\x)]&\leq \frac{2\alpha\beta}{(1-\beta)^2}(f(\x_0) - f(\x_{t})) + \|\x_0 - \x\|^2\\
& + \frac{\alpha^2}{(1-\beta)^2}(2s\beta+1)(G^2+\delta^2)(t+1)
\end{align*}
Let $\x=\x_*$ and noting that $f(\x_t)\geq f(\x_*)$, we have
\begin{align*}
\sum_{k=0}^t\E[f(\x_k) - f(\x_*)]\leq& \frac{\beta}{1-\beta}(f(\x_0) - f(\x_*)) \\
&+\frac{1-\beta}{2\alpha} \|\x_0 - \x_*\|^2 + \frac{\alpha}{2(1-\beta)}(2s\beta+1)(G^2+\delta^2)(t+1)
\end{align*}
Define $\xh_t = \sum_{k=0}^t\x_k/(t+1)$. By convexity of $f(\x)$, we have  
\[
\E[f(\xh_t) - f(\x_*)]\leq \frac{\beta}{(1-\beta)(t+1)}(f(\x_0) - f(\x_*)) + \frac{(1-\beta)\|\x_0 - \x_*\|^2}{2\alpha (t+1)} + \frac{\alpha(2s\beta+1)(G^2 +\delta^2)}{2(1-\beta)}
\]
By plugging the value of $\alpha$, we complete the proof of  Theorem~\ref{thm:1}.

\subsubsection{Proof of Theorem~\ref{thm:2}}
The proof of Theorem~\ref{thm:2} proceeds similarly as that of Theorem~\ref{thm:1} until we get~(\ref{eqn:later}), then we proceed as follows
\begin{align*}
&\|\x_{k+1} + \p_{k+1} - \x\|^2 = \|\x_k + \p_k - \x\|^2 - \frac{2\alpha}{1-\beta}(\x_k + \p_k - \x)^{\top}\G_k+ \left(\frac{\alpha}{1-\beta}\right)^2\|\G_k\|^2\\
& =  \|\x_k + \p_k - \x\|^2 - \frac{2\alpha}{1-\beta}(\x_k  - \x)^{\top}\G_k- \frac{2\alpha\beta}{(1-\beta)^2}(\x_k - \x_{k-1})^{\top}\G_k\\
& + \frac{2s\alpha^2\beta}{(1-\beta)^2}(\G_k - \G_{k-1})^{\top}\G_k -  \frac{2s\alpha^2\beta}{(1-\beta)^2}\|\G_k\|^2+ \left(\frac{\alpha}{1-\beta}\right)^2\|\G_k\|^2
\end{align*}
The above inequality holds for any $k\geq 0$ if we define $\G_{-1}=0$. 
Since $s\beta>1/2$ and $\|\G_k\|\leq M$, then 
\begin{align*}
\|\x_{k+1} + \p_{k+1} - \x\|^2& \leq   \|\x_k + \p_k - \x\|^2 - \frac{2\alpha}{1-\beta}(\x_k  - \x)^{\top}\G_k- \frac{2\alpha\beta}{(1-\beta)^2}(\x_k - \x_{k-1})^{\top}\G_k \\
&+ \frac{2s\alpha^2\beta}{(1-\beta)^2}M\|\G_k - \G_{k-1}\|
\end{align*}
Taking expectation on both sides and using the first two inequalities in~(\ref{eqn:fact}), we have
\begin{align*}
\sum_{k=0}^t\E[f(\x_k) - f(\x)]\leq \frac{\beta}{1-\beta}(f(\x_0) - f(\x)) +\frac{1-\beta}{2\alpha} \|\x_0 - \x\|^2 + \frac{s\alpha\beta M}{(1-\beta)}\sum_{t=0}^t\|\G_k - \G_{k-1}\|
\end{align*}
Then
\begin{align*}
\E[f(\xh_t) - f(\x)]\leq \frac{\beta(f(\x_0) - f(\x))}{(1-\beta)(t+1)} +\frac{(1-\beta) \|\x_0 - \x\|^2}{2\alpha(t+1)} + \frac{\alpha s\beta M}{(1-\beta)}\frac{\sum_{t=0}^t\|\G_k - \G_{k-1}\|}{t+1}
\end{align*}

\subsubsection{Proof of Theorem~\ref{thm:3}}
We first present two lemmas. 
\begin{lemma}\label{lem:1}
Let $\z_k = \x_k + \p_k$. For SUM,  we have for any $k\geq 0$, 
\begin{align*}
\E[f(\z_{k+1})- f(\z_k)]& \leq\frac{1}{2L}\E [\| \nabla f(\z_k) - \nabla f(\x_k)\|^2]\\
&+ \left(\frac{L\alpha^2}{(1-\beta)^2} -  \frac{\alpha}{1-\beta}\right)\E[\|\nabla f(\x_k)\|^2] + \frac{L\alpha^2}{2(1-\beta)^2}\sigma^2
\end{align*}
\end{lemma}
\begin{lemma}\label{lem:2}
For SUM, we have for any  $k\geq 0$, 
\begin{align*}
\E[\|\nabla f(\z_k) - \nabla f(\x_k)\|^2]\leq   \frac{L^2\beta^2((1-\beta)s -1)^2\alpha^2(G^2+\sigma^2)}{(1-\beta)^4}
\end{align*}
\end{lemma}

We continue the proof of Theorem~\ref{thm:3} as follows. Let $B, B'$ be defined as
\[
B =  \frac{\alpha}{(1-\beta)} - \frac{L\alpha^2}{(1-\beta)^2}, \quad B'=\frac{L\beta^{2}((1-\beta)s-1)^2\alpha^2(G^2+\sigma^2)}{2(1-\beta)^4} + \frac{L\alpha^2\sigma^2}{2(1-\beta)^2}
\]
Lemma~\ref{lem:1} and Lemma~\ref{lem:2} imply that
\begin{align*}
\E[f(\z_{k+1}) - f(\z_k)]\leq -B\E[\|\nabla f(\x_k)\|^2]  +B'
\end{align*}
By summing the above inequalities for $k=0,\ldots, t$ and noting that $\alpha< \frac{1-\beta}{L}$, 
\begin{align*}
B\sum_{k=0}^t\E[\|\nabla f(\x_k)\|^2]&\leq \E[f(\z_0) - f(\z_{t+1})] + (t+1)B'\leq \E[f(\z_0) - f_*] + (t+1)B'
\end{align*}
Then 
\begin{align*}
\min_{k=0,\ldots, t}\E[\|\nabla f(\x_k)\|^2]\leq \frac{f(\z_0) - f_*}{(t+1)B} + \frac{B'}{B}
\end{align*}
Assume $\alpha \leq \frac{1-\beta}{2L}$, then 
\[
B =\frac{\alpha}{1-\beta} - \frac{\alpha^2L }{(1-\beta)^2} = \frac{\alpha}{1-\beta}\left(1- \frac{\alpha L}{1-\beta}\right)\geq \frac{\alpha}{2(1-\beta)}
\]
Then 
\begin{align*}
\min_{k=0,\ldots, t}\E[\|\nabla f(\x_k)\|^2]\leq \frac{2(f(\z_0) - f_*)(1-\beta)}{\alpha (t+1)} + \frac{2(1-\beta)}{\alpha}B'
\end{align*}
Assume $\alpha =\min\{\frac{1-\beta}{2L}, \frac{C}{\sqrt{t+1}}\}$. Then
\begin{align*}
\min_{k=0,\ldots, t}\E[\|\nabla f(\x_k)\|^2]&\leq \frac{2(f(\z_0) - f_*)(1-\beta)}{t+1}\max\left\{\frac{2L}{1-\beta}, \frac{\sqrt{t+1}}{C}\right\}\\
& + \frac{C}{\sqrt{t+1}}\frac{L\beta^{2}((1-\beta)s-1)^2(G^2+\sigma^2) + L(1-\beta)^2\sigma^2}{(1-\beta)^3} 
\end{align*}
We then complete the proof by noting that $\z_0 = \x_0$. 

\subsubsection{Proof of Theorem~\ref{thm:4}}
With a slightly different analysis from that of Lemma~\ref{lem:1}, we can have the following lemma. 
\begin{lemma}\label{lem:3}
Let $\z_k = \x_k + \p_k$. For SUM,  we have for any $k\geq 0$, 
\begin{align*}
\E[f(\z_{k+1})- f(\z_k)]& \leq\frac{1}{2L((1-\beta)s-1)^2}\E [\| \nabla f(\z_k) - \nabla f(\x_k)\|^2]\\
&+ \left(\frac{[1+((1-\beta)s-1)^2]\alpha^2L}{2(1-\beta)^2} -  \frac{\alpha}{1-\beta}\right)\E[\|\nabla f(\x_k)\|^2] + \frac{L\alpha^2}{2(1-\beta)^2}\sigma^2
\end{align*}
\end{lemma}
With Lemma~\ref{lem:3} and Lemma~\ref{lem:2} and a similar analysis as that for Theorem~\ref{thm:3}, we can easily prove Theorem~\ref{thm:4}. 

\section{Empirical Studies}\label{sec:exp}
In this section, we present some empirical results. We focus on the non-convex optimization of  deep neural networks. 
We use two benchmark datasets, namely  CIFAR-10 and CIFAR-100,  and learn a deep convolutional neural network (CNN) for classification on the two datasets, respectively. Both datasets contain $50000$ training images of size $32*32$ from $10$ classes (CIFAR-10) or $100$ classes (CIFAR-100) and $10000$ testing images of the same size. The employed CNN consists of 3 convolutional layers and 2 fully-connected layers. Each convolutional layer is followed by a max pooling layer. The output of the last fully-connected layer is fed into a $10$-class or $100$-class softmax loss function.  We emphasize that we do not intend to obtain the state-of-the-art prediction performance by trying different network structures and different engineering tricks, but instead focus our attention on optimization.  We compare the three variants of SUM, i.e., SHB, SNAG, and SG, which corresponds to $s=0$, $s=1$ and $s = 1/(1-\beta)$ in~(\ref{eqn:um}). We fix the momentum constant $\beta=0.9$ and the regularization parameter of weights to $0.0005$. We use a mini-batch of size $128$ to compute a stochastic gradient at each iteration.  All three methods use the same initialization.  We follow the procedure in~\citep{krizhevsky2012imagenet} to set the step size $\alpha$, i.e., initially giving a relatively large step size and and decreasing the step size by $10$ times after certain number of iterations when observing the performance on testing data saturates.

\begin{figure}[t]
\begin{center}
\includegraphics[scale=0.25]{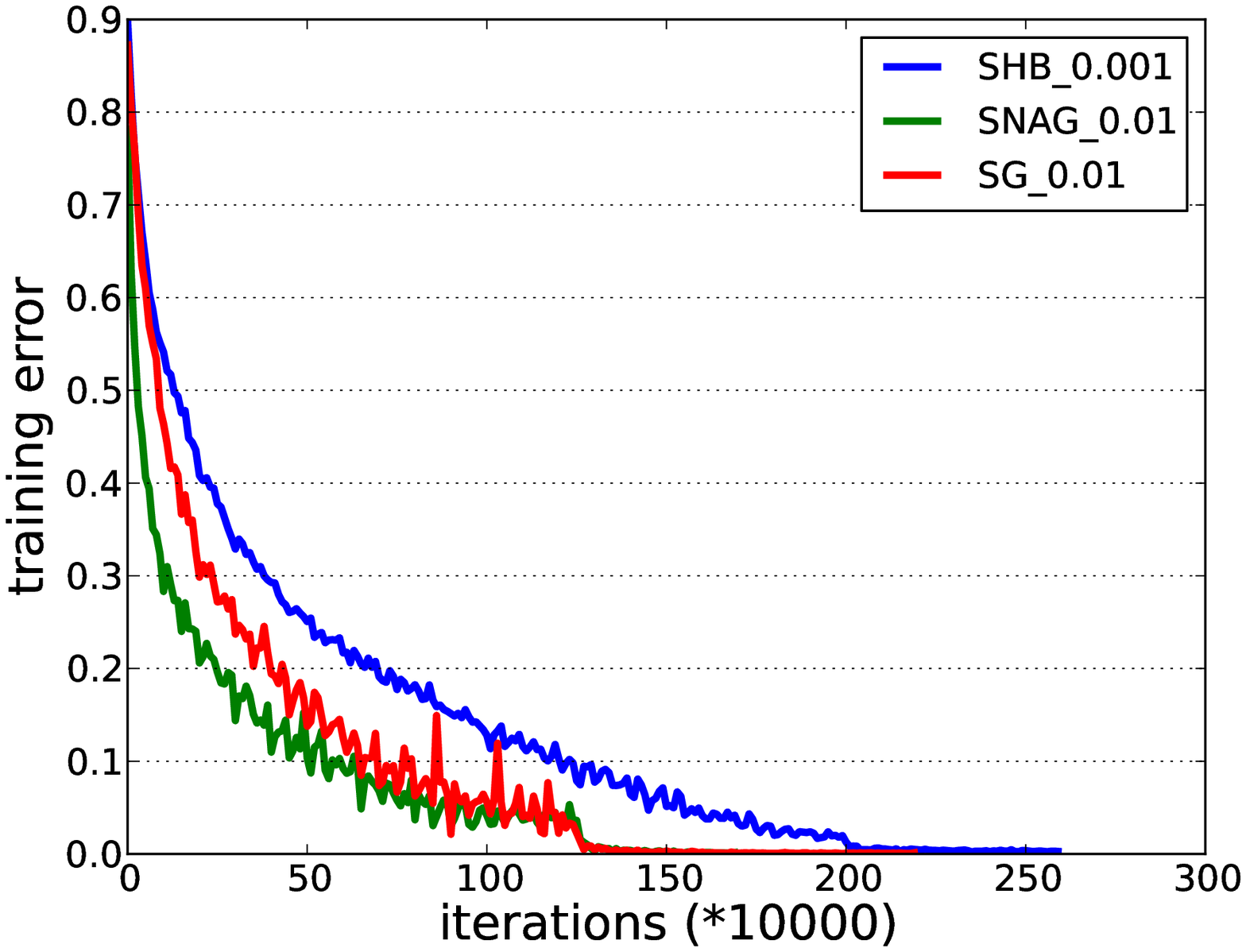}
\includegraphics[scale=0.25]{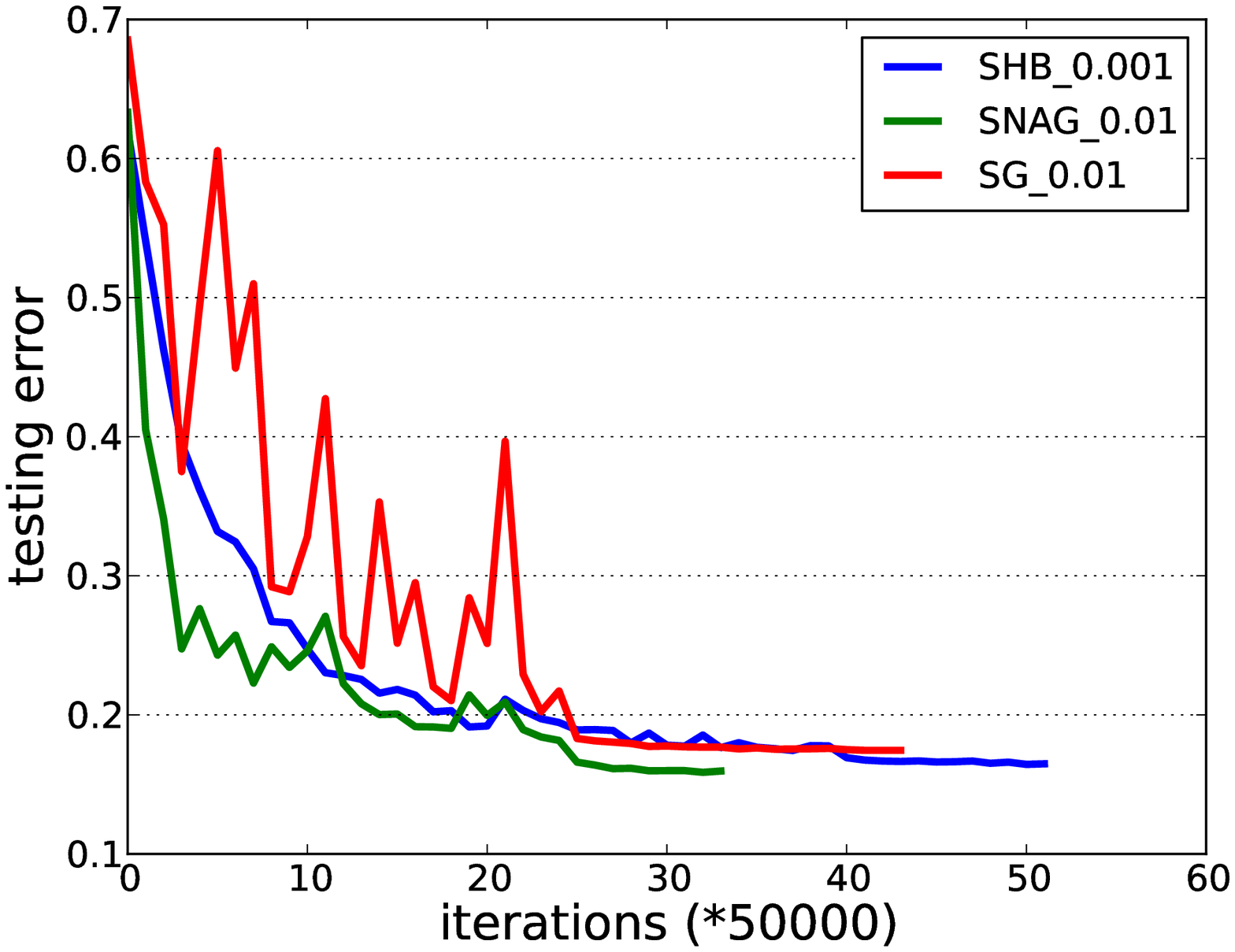}

\caption{Training and testing error of different methods with the best initial step size on CIFAR-10.}
\label{fig:1}
\end{center}

\begin{center}
\includegraphics[scale=0.25]{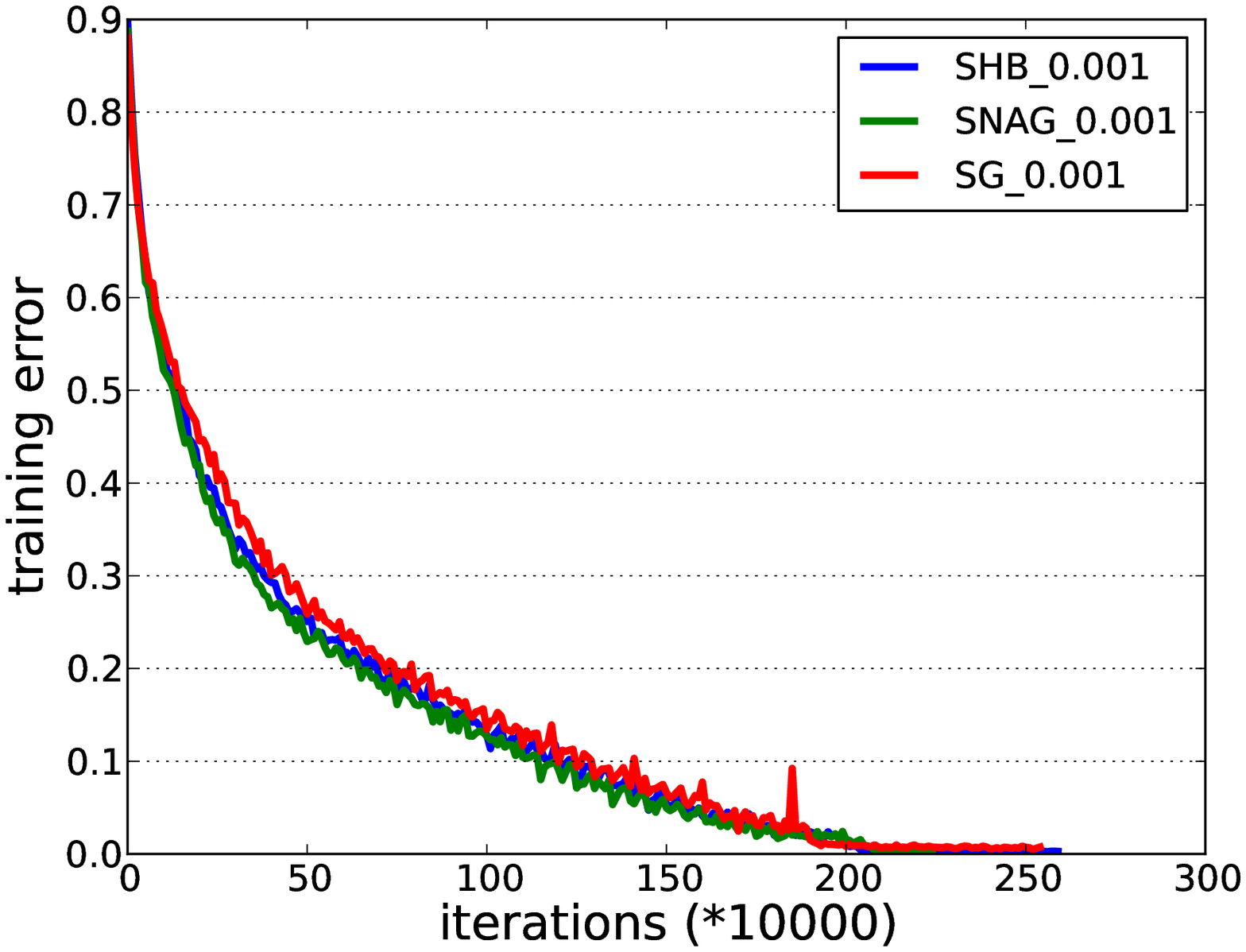}
\includegraphics[scale=0.25]{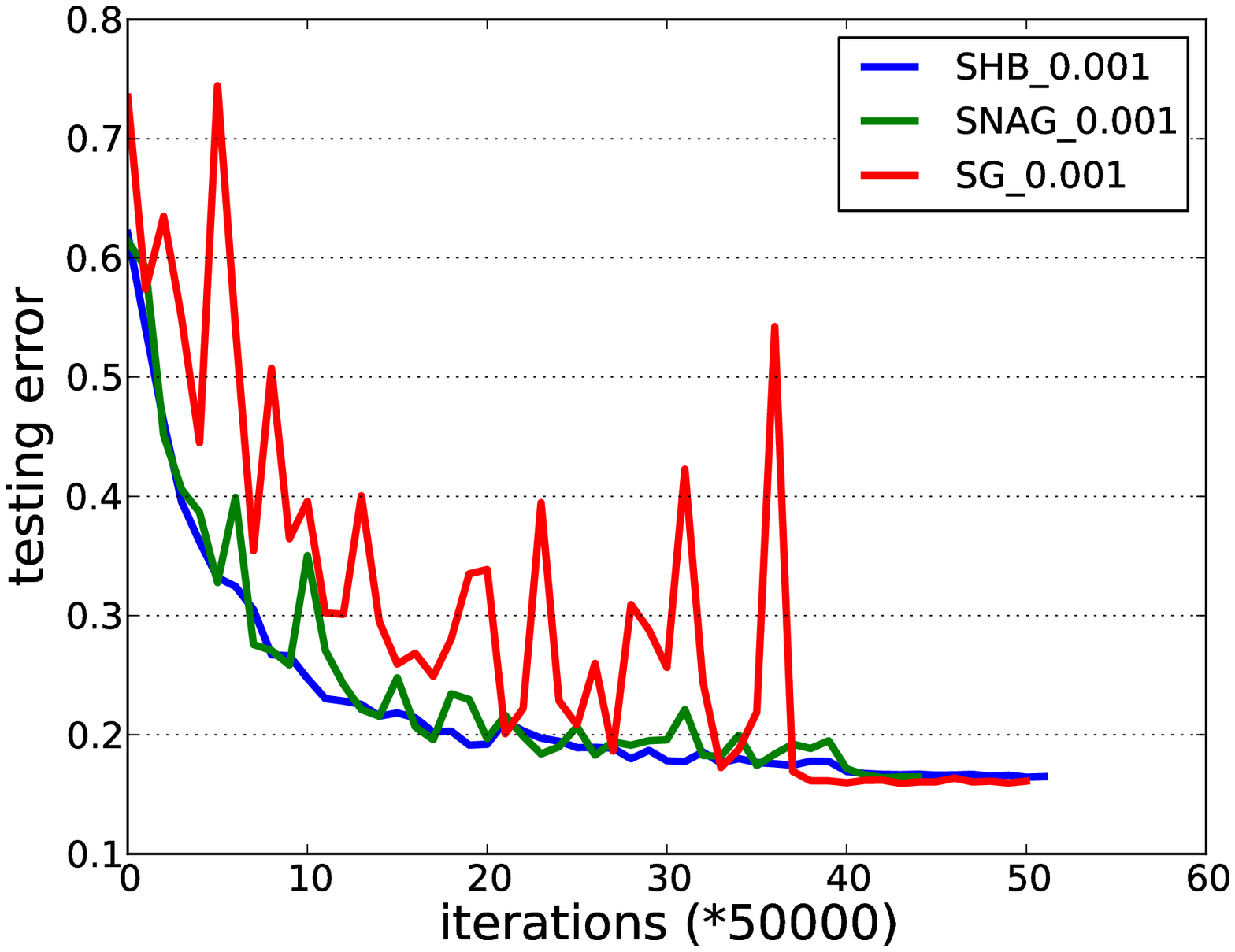}
\caption{Training and testing error of different methods with the same initial step size on CIFAR-10.}
\label{fig:2}
\end{center}
\end{figure}

First, we report some experimental results on CIFAR-10 data. For the initial step size, we search in a range ($\{0.001, 0.002, 0.005, 0.01, 0.02, 0.05\}$) for different methods and select the best one that yields the fastest convergence in training error. In particular,  for SHB the best initial step size is $0.001$ and that for SNAG and SG is $0.01$. In fact, using a larger initial step size (e.g, $0.002$) for SHB gives a divergent result. The training and testing error of different methods versus the number of iterations is plotted in Figure~\ref{fig:1}. We also plot the performance of different methods with the same initial step size $0.001$ in Figure~\ref{fig:2}. From the results, we can see that with the same step size $0.001$ the three methods converge similarly in training error, which is consistent with the result in Theorem~\ref{thm:3}. However, the testing error curves of different methods exhibit different degree of oscillation. In particular, the degree of oscillation of testing error across iterations for the three methods is ordered in SG, SNAG, SHB from large to small, which coincides with our unified framework for explaining the three methods with the free parameter $s$ varying from $10$~\footnote{which is the value of $1/(1-\beta)$ because the momentum constant is $\beta=0.9$.} to $1$ and $0$. After the step size is decreased to a very small value that can discount the large value of $s$, the oscillation of testing error for  all methods is gone~\footnote{See SG after $25$ epochs in Figure~\ref{fig:1}.}.  Moreover, we observe that SG and SNAG allow a larger initial step size ($0.01$) on this data, which makes them converge faster in training error without obvious loss in testing performance~\footnote{Indeed SNAG achieves better testing error.}. This result also verifies  Theorem~\ref{thm:4}. 

\begin{figure}[t]
\centering
\includegraphics[scale=0.25]{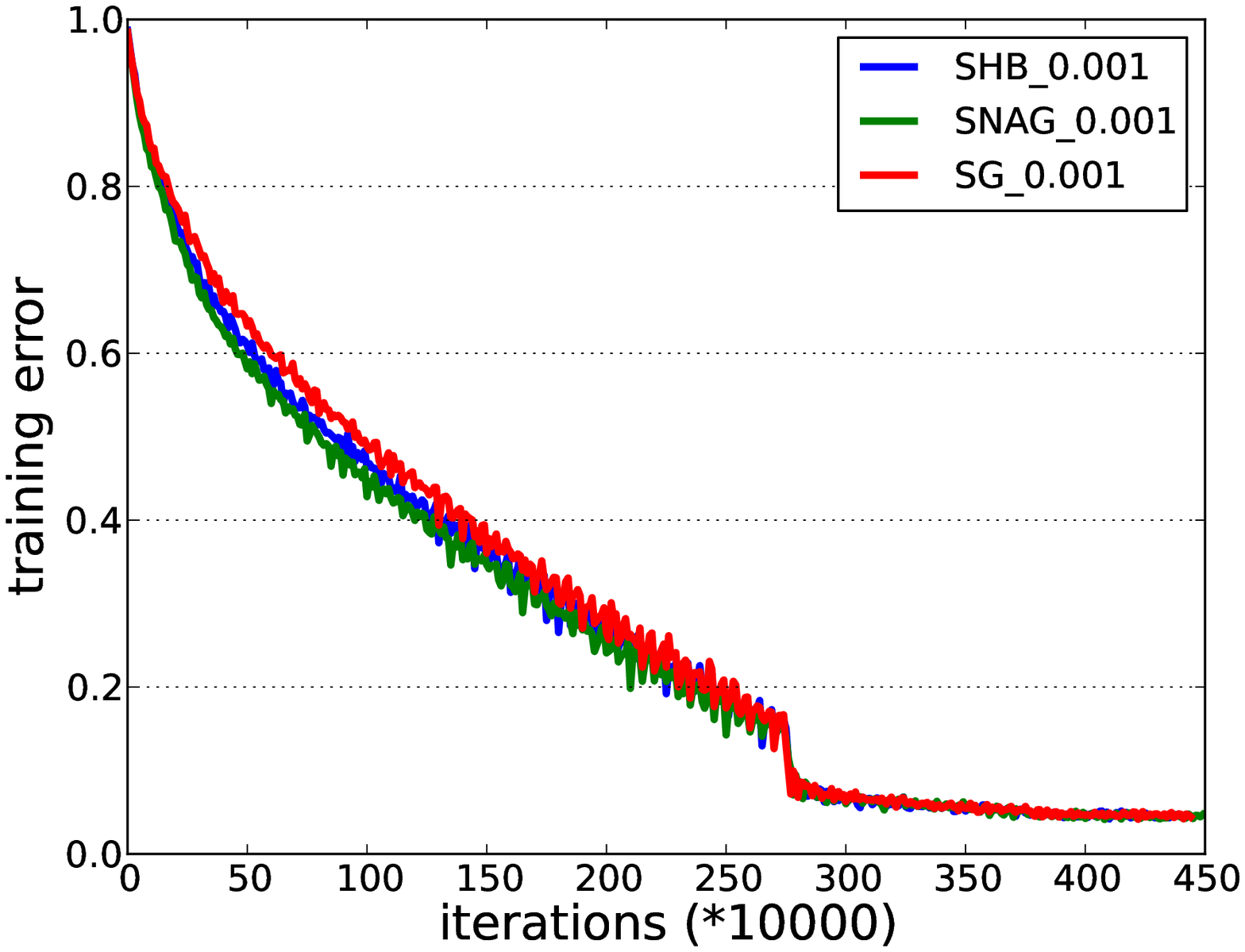}
\includegraphics[scale=0.25]{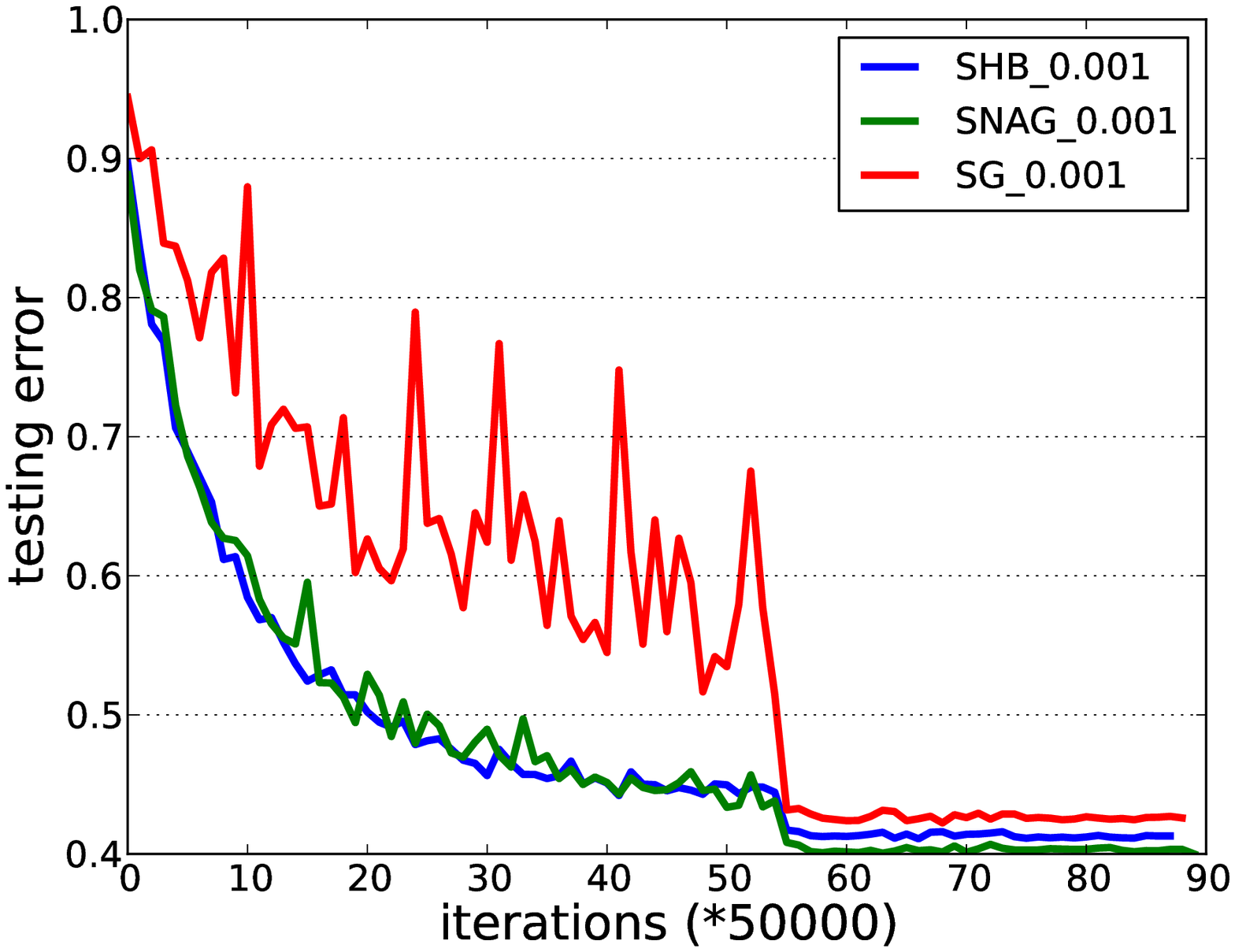}
\caption{Training and testing error of different methods with the initial step size $0.001$ on CIFAR-100.}
\label{fig:3}
\includegraphics[scale=0.25]{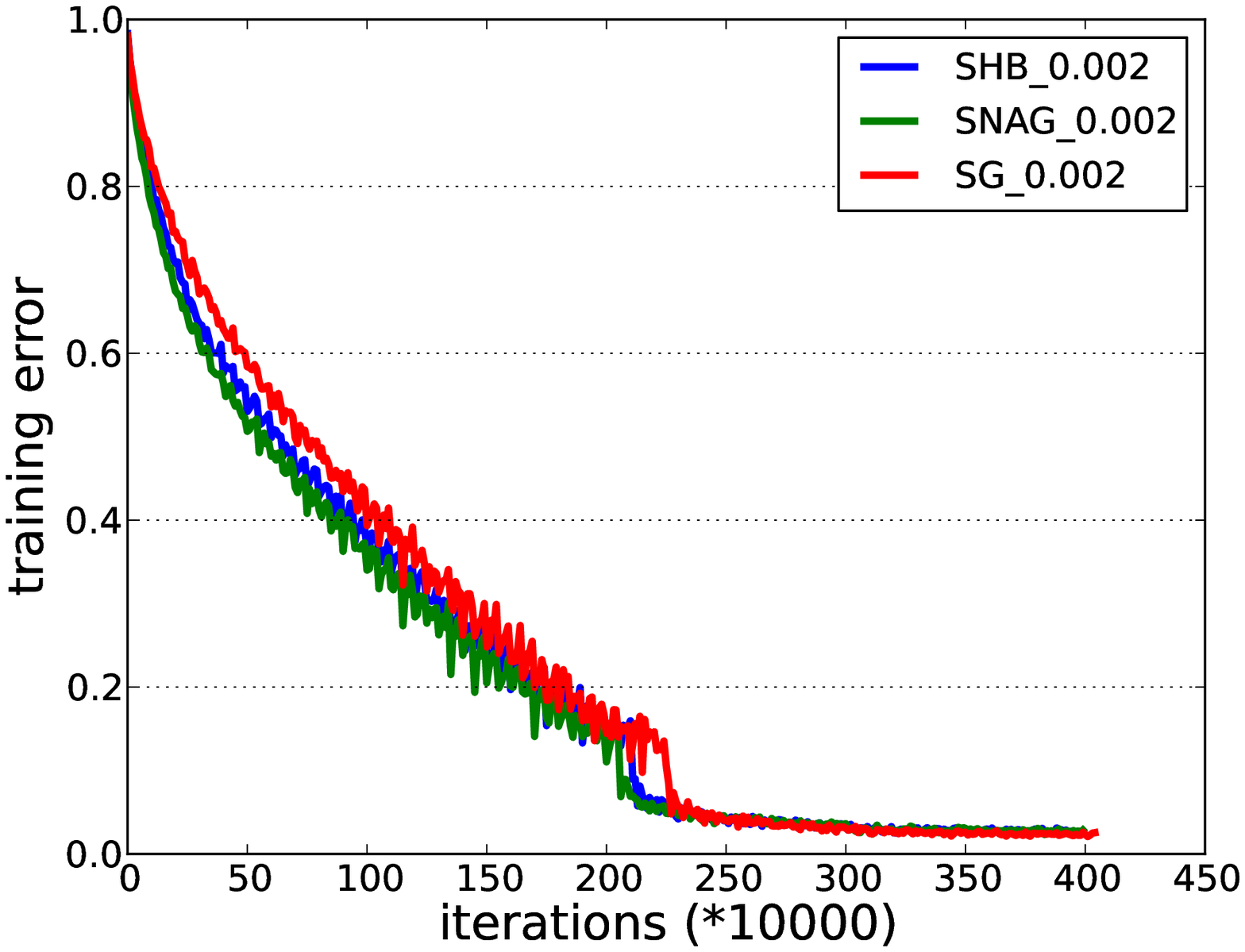}
\includegraphics[scale=0.25]{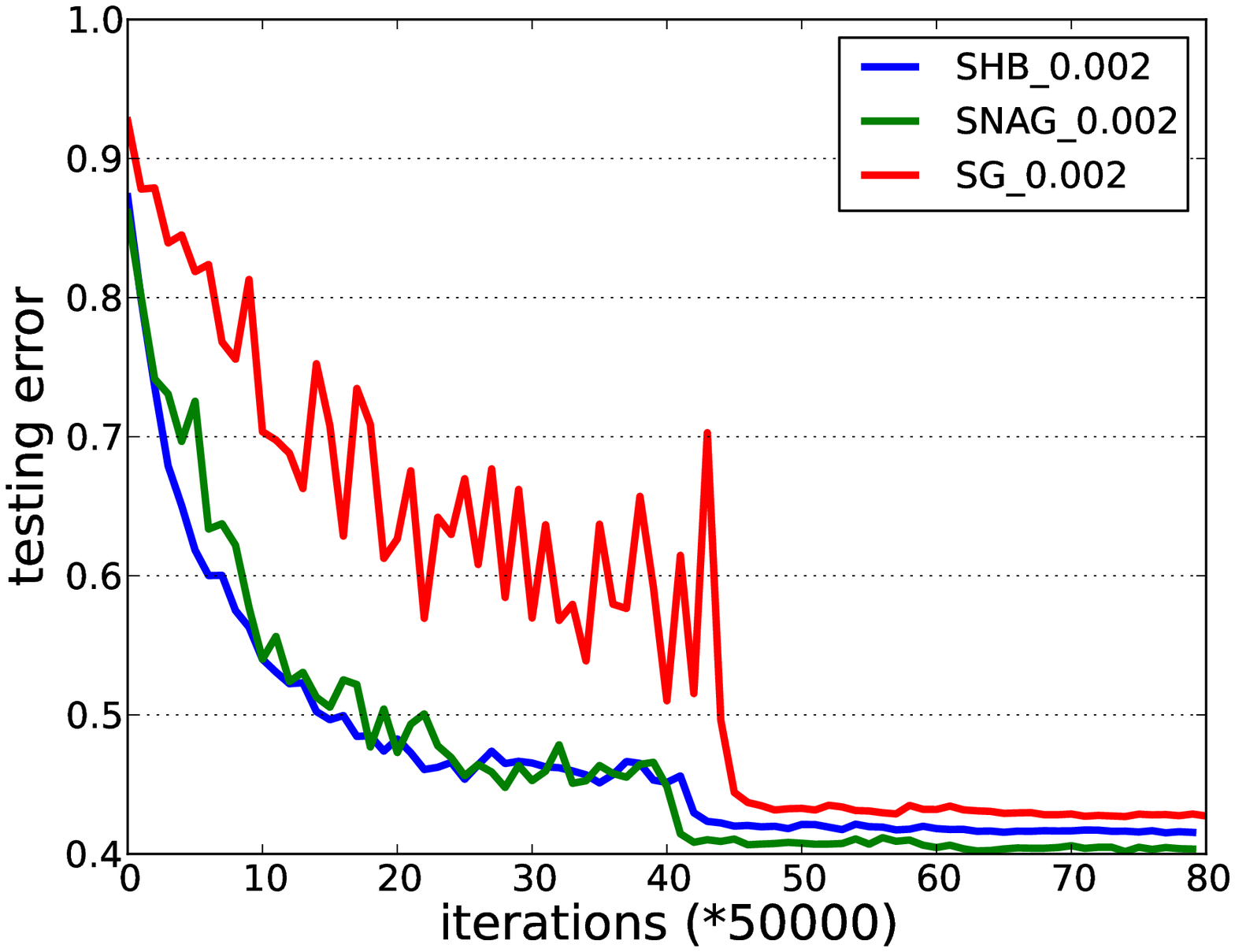}
\caption{Training and testing error of different methods with the initial step size $0.002$ on CIFAR-100.}
\label{fig:4}
\end{figure}

\begin{figure}[t]
\centering
\includegraphics[scale=0.25]{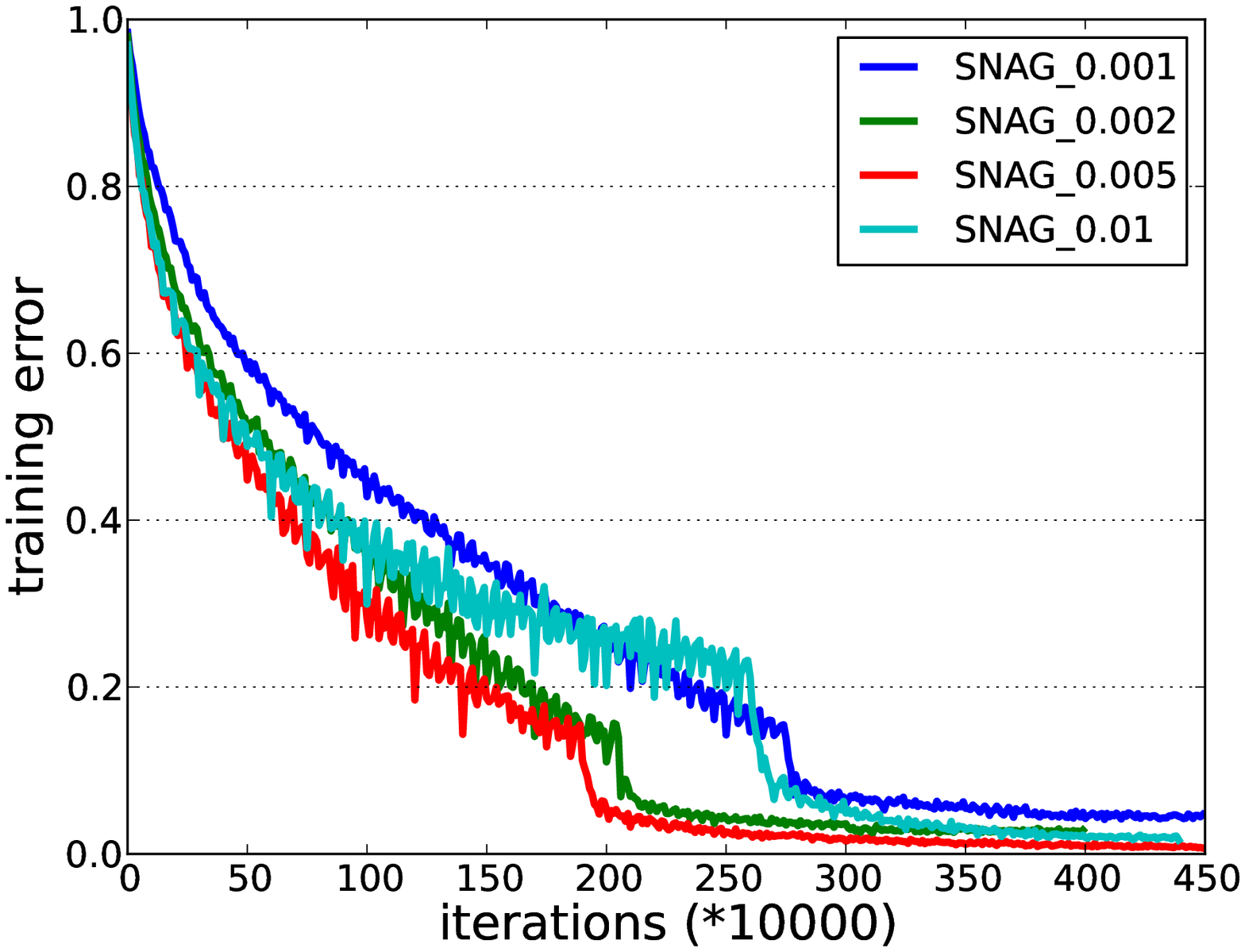}
\includegraphics[scale=0.25]{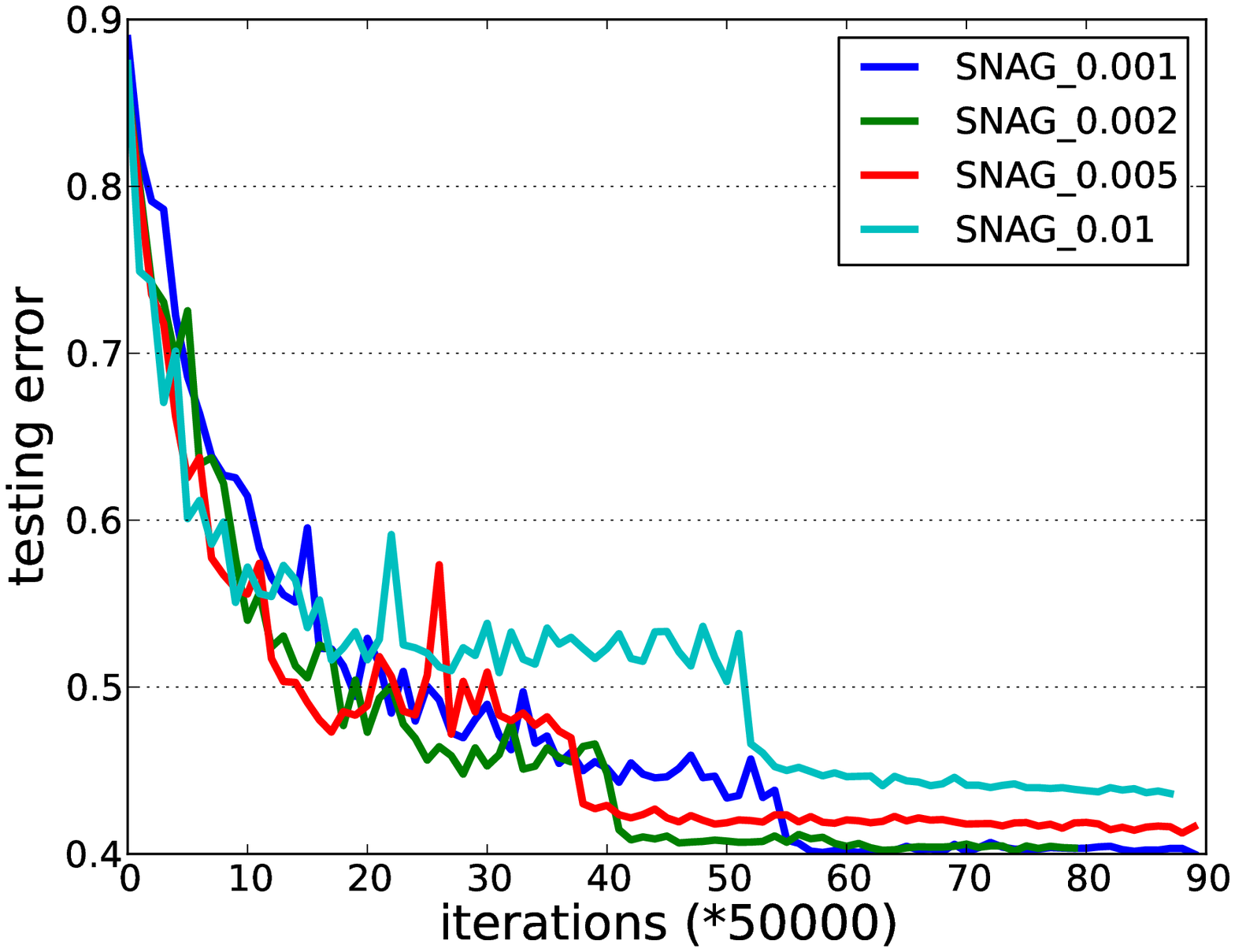}
\caption{Training and testing error of SNAG with different  initial step size on CIFAR-100.}
\label{fig:5}
\end{figure}

Second, we report some experimental results on CIFAR-100 data.  We plot the training and testing error of the three methods with two different initial step sizes in Figure~\ref{fig:3} and Figure~\ref{fig:4}, respectively.  We observe similar trends in the training error and testing error, i.e., the three methods have similar performance (convergence speed) in training error but exhibits different degree of oscillation in the testing error convergence. They also show that SNAG achieves the best prediction performance on the testing data. For this data, we also  observe that a larger step size $0.005$ yields a faster convergence in the training error for SNAG, but it gives worse  testing error, which is illustrated in Figure~\ref{fig:5}.  

\begin{figure}[t]
\centering
\includegraphics[scale=0.25]{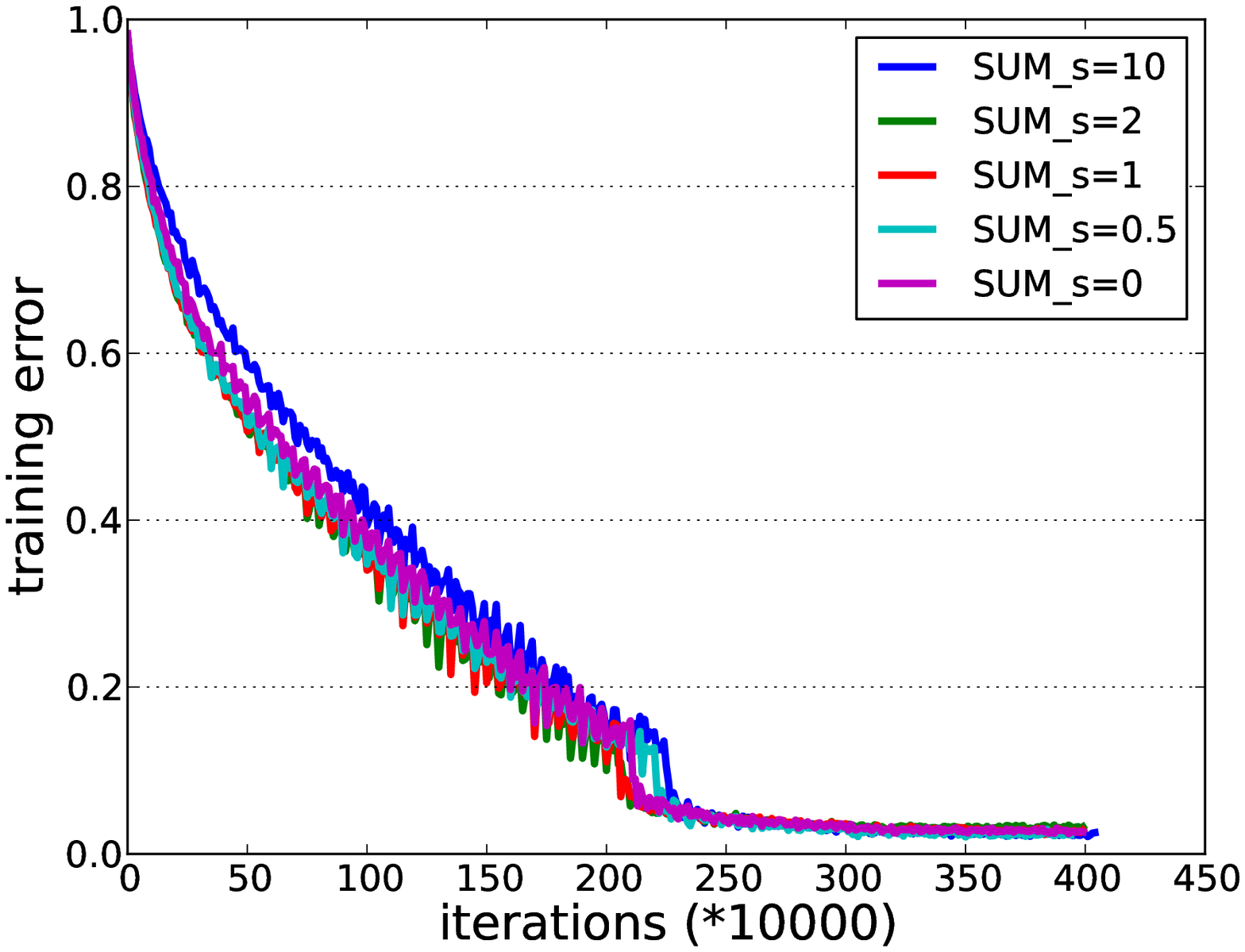}
\includegraphics[scale=0.25]{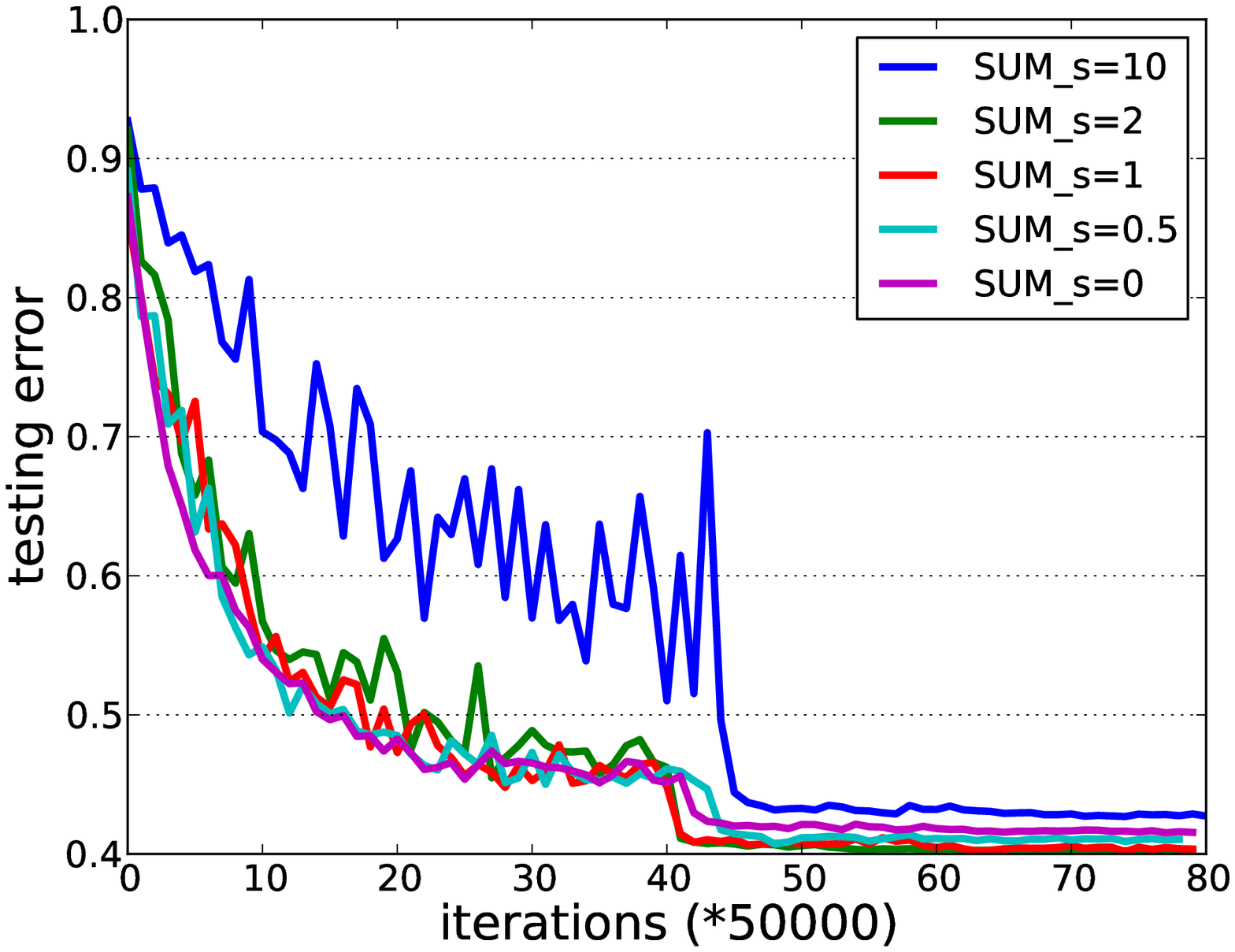}
\caption{Training and testing error of SUM with different  $s$ on CIFAR-100.}
\label{fig:6}
\end{figure}

Finally, we present a comparison of SUM with different values of $s$ including variants besides SHB, SNAG  and SG. In particular, we compare SUM with $s\in\{0, 0.5, 1, 2, 10\}$. Note that $s=0$ corresponds to SHB, $s=1$ corresponds to SNAG and $s=10$ corresponds to SG since $\beta=0.9$. The results on the CIFAR-100 data are shown in Figure~\ref{fig:6}. From the results, we can observe that for the convergence of  training error  the different variants perform similarly with SG ($s=10$) performing slightly worse. For the convergence of testing error, we observe a clear change from $s=10$ to $s=0$ in that the oscillation decreases.  However, SUM with  $s=2$ converge to a slightly better solution in terms of the testing error.

\section{Conclusion}
In this paper, we have developed a unified framework of momentum methods that subsumes the heavy-ball method, the Nesterov's accelerated gradient method, the gradient method and their stochastic variants as special cases.  We have also established basic convergence results of the stochastic momentum methods for both convex and non-convex optimization, which have been widely adopted in training deep neural networks. The unified framework exhibits that the difference between the heavy-ball, the accelerated gradient method and the gradient method  lies at  the momentum term between two consecutive auxiliary points (i.e., the step size used to compute the auxiliary sequence). The unified framework helps explain some experimental results for optimization deep neural networks, but also raises some open problems, e.g.,   how to quantitatively analyze the theoretical impact of different  momentum terms on the convergence of generalization error for deep learning?

\bibliographystyle{plain}

 \bibliography{all}
 
 \appendix
 \section*{Proof of Lemma~\ref{lem:0}}
\begin{proof}
For clarity of presentation, we let $\E_k[\cdot]$ denote the expectation over the randomness in $\xi_1,\ldots, \xi_k$ and let $\E_{k|k-1}[\cdot]$ denote the expectation over $\xi_k$ given $\xi_1,\ldots, \xi_{k-1}$ fixed. We use $\E[\cdot]$ to denote the expectation over all randomness. We abuse the notation $\partial f(\x_k)$ to denote $\E_{k|k-1}[\G_k]$. Define $\delta_k =\G_k- \partial f(\x_k)$. 
From the recursive relation in~(\ref{eqn:rec}) and the definition of $\p_k$ in~(\ref{eqn:p}), we have for any $\x\in\R^d$ and $k\geq 1$
\begin{align}\label{eqn:later}
&\|\x_{k+1} + \p_{k+1} - \x\|^2 = \|\x_k + \p_k - \x\|^2 - \frac{2\alpha}{1-\beta}(\x_k + \p_k - \x)^{\top}\G_k+ \left(\frac{\alpha}{1-\beta}\right)^2\|\G_k\|^2\nonumber\\
& =  \|\x_k + \p_k - \x\|^2 - \frac{2\alpha}{1-\beta}(\x_k  - \x)^{\top}\G_k- \frac{2\alpha\beta}{(1-\beta)^2}(\x_k - \x_{k-1})^{\top}\G_k\nonumber\\
& - \frac{2s\alpha^2\beta}{(1-\beta)^2}\G_{k-1}^{\top}\G_k+ \left(\frac{\alpha}{1-\beta}\right)^2\|\G_k\|^2\\
& =  \|\x_k + \p_k - \x\|^2 - \frac{2\alpha}{1-\beta}(\x_k  - \x)^{\top}(\delta_k+\partial f(\x_k)) - \frac{2\alpha\beta}{(1-\beta)^2}(\x_k - \x_{k-1})^{\top}(\delta_k+\partial f(\x_k))\nonumber\\
&- \frac{2s\alpha^2\beta}{(1-\beta)^2}(\delta_{k-1}+\partial f(\x_{k-1}))^{\top}(\delta_k + \partial f(\x_k))+ \left(\frac{\alpha}{1-\beta}\right)^2\|\delta_k + \partial f(\x_k)\|^2\nonumber 
\end{align}
Taking expectation over both sides and by noting that 
\begin{align*}
&\E_k[(\x_k  - \x)^{\top}(\delta_k+\partial f(\x_k)) ] = \E_{k-1}[(\x_k  - \x)^{\top}\partial f(\x_k)] \\
&\E_k[(\x_k - \x_{k-1})^{\top}(\delta_k+\partial f(\x_k))] = \E_{k-1}[(\x_k - \x_{k-1})^{\top}\partial f(\x_k)]\\
&\E_k[(\delta_{k-1}+\partial f(\x_{k-1}))^{\top}(\delta_k + \partial f(\x_k))] = \E_k[\delta_{k-1}^{\top}\delta_k]+ \E_k[\delta_{k-1}^{\top}\partial f(\x_k)] + \E_k[\delta_{k}^{\top}\partial f(\x_{k-1})]\\
&\hspace*{2in} + \E_k[\partial f(\x_{k-1})^{\top}\partial f(\x_k)] = \E_{k-1}[[(\delta_{k-1}+\partial f(\x_{k-1}))^{\top}\partial f(\x_k)]\\
&\E_k[\|\delta_k + \partial f(\x_k)\|^2] = \E_k[\|\delta_k\|^2] + \E_{k-1}[\|\partial f(\x_k)\|^2]
\end{align*}
we have
\begin{align*}
&\E_k[\|\x_{k+1} + \p_{k+1} - \x\|^2] = \E_{k-1}[\|\x_k + \p_k - \x\|^2] \\
&- \frac{2\alpha}{1-\beta}\E_{k-1}[(\x_k  - \x)^{\top}\partial f(\x_k)] - \frac{2\alpha\beta}{(1-\beta)^2}\E_{k-1}[(\x_k - \x_{k-1})^{\top}\partial f(\x_k)]\\
&- \frac{2s\alpha^2\beta}{(1-\beta)^2}\E_{k-1}[\G_{k-1}^{\top}\partial f(\x_k)]+ \left(\frac{\alpha}{1-\beta}\right)^2(\E_k[\|\delta_k\|^2 + \E_{k-1}[\|\partial f(\x_k)\|^2 ])
\end{align*}
We also note that 
\begin{equation}\label{eqn:fact}
\begin{aligned}
&f(\x_k ) - f(\x)\leq (\x_k - \x)^{\top}\partial f(\x_k)\\
&f(\x_k) - f(\x_{k-1})\leq (\x_k - \x_{k-1})^{\top}\partial f(\x_k)\\
&-\E[\G_{k-1}^{\top}\partial f(\x_k)]\leq \frac{\E[\|\G_{k-1}\|^2 + \|\partial f(\x_k)\|^2]}{2}\leq \delta^2/2 + G^2\leq \delta^2 + G^2\\
&\E_k[\|\delta_k\|^2]\leq \delta^2,\quad\E_{k-1}[\|\partial f(\x_k)\|^2\leq G^2
\end{aligned}
\end{equation}
where the first two inequalities are due to the convexity and the last three inequalities are due to the boundness assumption.  Thus
\begin{align*}
\E[\|\x_{k+1} + \p_{k+1} - \x\|^2&] \leq   \E[\|\x_k + \p_k - \x\|^2] - \frac{2\alpha}{1-\beta}\E[(f(\x_k) - f(\x))]\\
 &- \frac{2\alpha\beta}{(1-\beta)^2}\E[(f(\x_k) - f(\x_{k-1}))]  + \left(\frac{\alpha}{1-\beta}\right)^2(2s\beta+1) (G^2 +\delta^2)\quad k\geq 1
\end{align*}
By letting  $\x_{-1}=\x_0$, we note that the above inequality also holds for $k=0$ following a similar analysis.
\end{proof}
 
 \section*{Proof of Lemma~\ref{lem:1}}
\begin{proof}
Since $f(\x)$ is a smooth function, we have
\begin{align*}
 f(\y) - f(\x) - \nabla f(\x)^{\top}(\y - \x)\leq \frac{L}{2}\|\y - \x\|^2
\end{align*}
Then
\begin{align*}
f(\z_{k+1})& \leq f(\z_k) + \nabla f(\z_k)^{\top}(\z_{k+1} - \z_k) + \frac{L}{2}\frac{\alpha^2}{(1-\beta)^2}\|\G_k\|^2\\
& = f(\z_k) - \frac{\alpha}{1-\beta}\nabla f(\z_k)^{\top}\G_k +  \frac{L}{2}\frac{\alpha^2}{(1-\beta)^2}\|\G_k\|^2\\
&  = f(\z_k) - \frac{\alpha}{1-\beta}\nabla f(\z_k)^{\top}(\delta_k+\nabla f(\x_k)) +  \frac{L}{2}\frac{\alpha^2}{(1-\beta)^2}\|\G_k\|^2\\
&  = f(\z_k) - \frac{\alpha}{1-\beta}\nabla f(\z_k)^{\top}\delta_k - \frac{\alpha}{1-\beta} \nabla f(\z_k)^{\top}\nabla f(\x_k)+ \frac{L}{2}\frac{\alpha^2}{(1-\beta)^2}\|\G_k\|^2\\
&  = f(\z_k) - \frac{\alpha}{1-\beta}\nabla f(\z_k)^{\top}\delta_k - \frac{\alpha}{1-\beta} (\nabla f(\z_k) - \nabla f(\x_k))^{\top}\nabla f(\x_k) - \frac{\alpha}{1-\beta}\|\nabla f(\x_k)\|^2 \\
&+ \frac{L}{2}\frac{\alpha^2}{(1-\beta)^2}\|\G_k\|^2
\end{align*}

Taking expectation on both sides 
\begin{align}
&\E[f(\z_{k+1})- f(\z_k)]\leq \E\left[ -\frac{\alpha}{1-\beta} (\nabla f(\z_k) - \nabla f(\x_k))^{\top}\nabla f(\x_k) -  \frac{\alpha}{1-\beta}\|\nabla f(\x_k)\|^2  +\frac{L}{2}\frac{\alpha^2}{(1-\beta)^2}\|\nabla f(\x_k)\|^2 \right]\nonumber\\
&  +\frac{L}{2}\frac{\alpha^2}{(1-\beta)^2}\E[\|\G_k- \nabla f(\x_k)\|^2]\nonumber\\
&= \E\left[ -\frac{\alpha}{1-\beta} (\nabla f(\z_k) - \nabla f(\x_k))^{\top}\nabla f(\x_k) \right] + \left(\frac{L}{2}\frac{\alpha^2}{(1-\beta)^2} -  \frac{\alpha}{1-\beta}\right)\E[\|\nabla f(\x_k)\|^2] + \frac{L\alpha^2}{2(1-\beta)^2}\sigma^2\\
&\leq\frac{1}{2}\E \left[\frac{1}{L}\| \nabla f(\z_k) - \nabla f(\x_k)\|^2 + \frac{L\alpha^2}{(1-\beta)^2}\|\nabla f(\x_k)\|^2\right]\nonumber\\
&+ \left(\frac{L}{2}\frac{\alpha^2}{(1-\beta)^2} -  \frac{\alpha}{1-\beta}\right)\E[\|\nabla f(\x_k)\|^2] + \frac{L\alpha^2}{2(1-\beta)^2}\sigma^2\nonumber
\end{align}
\end{proof}

\section*{Proof of Lemma~\ref{lem:2}}
\begin{proof}
%
\begin{align*}
\|\nabla f(\z_k) - \nabla f(\x_k)\|^2\leq L^2\|\z_k - \x_k\|^2 = L^2 \|\p_k\|^2 =  \frac{L^2\beta^2}{(1-\beta)^2}\left\|\x_k  - \x_{k-1} + s\alpha\G_{k-1}\right\|^2
\end{align*}
Recall the recursion in~(\ref{eqn:rec2}). 
\[
\v_{k+1} = \beta\v_k +  ((1-\beta)s-1)\alpha \G_k
\]
Note that $\v_0=0$.  Denote by $\hat\alpha = \alpha ((1-\beta)s-1)$. 
By induction,  we can show that
\[
\v_k = \hat\alpha\sum_{i=0}^{k-1}\beta^{k-1-i}\G_i=  \hat\alpha \sum_{i=0}^{k-1}\beta^{i}\G_{k-1-i}
\]
Let $\Gamma_{k-1}= \sum_{i=0}^{k-1}\beta^i = \frac{1-\beta^k}{1-\beta}$. Then
\begin{align*}
\|\v_k\|^2&=\left\|\sum_{i=0}^{k-1}\frac{\beta^i}{\Gamma_{k-1}}\hat\alpha \G_{k-1-i}\right\|^2\Gamma_{k-1}^2\leq\Gamma_{k-1}^2 \sum_{i=0}^{k-1}\frac{\beta^i}{\Gamma_{k-1}}\hat\alpha^2\|\G_{k-1-i}\|^2=\Gamma_{k-1}\sum_{i=0}^{k-1}\beta^i\hat\alpha^2\|\G_{k-1-i}\|^2
\end{align*}
Then
\begin{align*}
\E[\|\v_k\|^2]\leq \Gamma_{k-1}\sum_{i=0}^{k-1}\beta^i\hat\alpha^2(G^2 + \sigma^2)= \Gamma_{k-1}^2\hat\alpha^2(G^2+\sigma^2)\leq\frac{\alpha^2((1-\beta)s-1)^2(G^2+\sigma^2)}{(1-\beta)^2}
\end{align*}
Then 
\begin{align*}
\|\nabla f(\z_k) - \nabla f(\x_k)\|^2&\leq \frac{L^2\beta^2}{(1-\beta)^2}\E[\|\x_k - \x_{k-1} + s\alpha \G_{k-1}\|^2] = \frac{L^2\beta^2}{(1-\beta)^2}\E[\|\v_k\|^2]\\
& \leq \frac{L^2\beta^2((1-\beta)s-1)^2\alpha^2(G^2+\sigma^2)}{(1-\beta)^4}
\end{align*}
\end{proof}

 \section*{Proof of Lemma~\ref{lem:3}}
\begin{proof}
We can follow the same analysis as in the proof of Lemma~\ref{lem:1} and get \begin{align*}
&\E[f(\z_{k+1})- f(\z_k)]\leq  \E\left[ -\frac{\alpha}{1-\beta} (\nabla f(\z_k) - \nabla f(\x_k))^{\top}\nabla f(\x_k) \right] \\
&+ \left(\frac{L}{2}\frac{\alpha^2}{(1-\beta)^2} -  \frac{\alpha}{1-\beta}\right)\E[\|\nabla f(\x_k)\|^2] + \frac{L\alpha^2}{2(1-\beta)^2}\sigma^2\\
&\leq\frac{1}{2}\E \left[\frac{1}{L((1-\beta)s-1)^2}\| \nabla f(\z_k) - \nabla f(\x_k)\|^2 + \frac{L\alpha^2((1-\beta)s-1)^2}{(1-\beta)^2}\|\nabla f(\x_k)\|^2\right]\\
&+ \left(\frac{L}{2}\frac{\alpha^2}{(1-\beta)^2} -  \frac{\alpha}{1-\beta}\right)\E[\|\nabla f(\x_k)\|^2] + \frac{L\alpha^2}{2(1-\beta)^2}\sigma^2\\
&=\frac{1}{2L((1-\beta)s-1)^2}\E[\| \nabla f(\z_k) - \nabla f(\x_k)\|^2]\\
&+ \left(\frac{\alpha^2L[1+((1-\beta)s-1)^2]}{2(1-\beta)^2} -  \frac{\alpha}{1-\beta}\right)\E[\|\nabla f(\x_k)\|^2] + \frac{L\alpha^2}{2(1-\beta)^2}\sigma^2
\end{align*}
\end{proof}

\end{document}